\documentclass[a4paper,11pt]{article}

\usepackage{natbib}
\usepackage{amsmath,amssymb,amsfonts,stmaryrd}
\usepackage[latin1]{inputenc}
\usepackage{color,rotating}
\usepackage[margin=25mm]{geometry}
\usepackage{authblk}
\usepackage{tikz}

\newcommand{\indep}{\,\rotatebox{90}{$\hspace{-.1cm}\models$}\,}
\newcommand{\matr}[1]{\boldsymbol{#1}}
\newcommand{\ud}{\mathrm{d}}
\newcommand{\vect}[1]{\boldsymbol{#1}}

\newcommand{\gcite}{\citep}
\newcommand{\gcitet}{\citet}

\newtheorem{theorem}{Theorem}

\graphicspath{{Figures/}}

\title{Exponential decay of pairwise correlation in Gaussian graphical models with an equicorrelational one-dimensional connection pattern}

\author{Guillaume Marrelec\textsuperscript{1,2,*}, Alain Giron\textsuperscript{1,2,*} and Laura Messio\textsuperscript{3,4,$\dagger$}\\
\footnotesize \textsuperscript{1} Sorbonne Universit{\'e}, CNRS, INSERM, Laboratoire d'imagerie biom{\'e}dicale, LIB, F-75006, Paris, France\\
\footnotesize \textsuperscript{2} Center for Interaction Science (CIS), Centre de recherches et d'études en sciences des interactions (Cr{\'e}si), F-75006, Paris, France\\
\footnotesize \textsuperscript{3} Sorbonne Université, CNRS, Laboratoire de Physique Théorique de la Matière Condensée, LPTMC, F-75005 Paris, France\\
\footnotesize \textsuperscript{4} Institut Universitaire de France, IUF, F-75005 Paris, France\\
\footnotesize \textsuperscript{*} Email: firstname.lastname@inserm.fr\\
\footnotesize \textsuperscript{$\dagger$} Email: laura.messio@lptmc.jussieu.fr}

\date{}

\begin{document}

\maketitle

\begin{abstract}
 We consider Gaussian graphical models associated with an equicorrelational and one-dimensional conditional independence graph. We show that pairwise correlation decays exponentially as a function of distance. We also provide a limit when the number of variables tend to infinity and quantify the difference between the finite and infinite cases.
 \par
 \
\end{abstract}

\section{Introduction}

Let $\vect{X} = ( X_1, \dots, X_n )$ be an $n$-dimensional variable. A conditional independence graph on $\vect{X}$ is a graphical representation of $\vect{X}$ which emphasizes the relationships of conditional independence between the $X_i$'s \gcite{Whittaker-1990}. More precisely, there is \emph{no} link between nodes $i$ and $j$ if $X_i$ and $X_j$ are conditionally independent given $\vect{X}_{[n] \setminus \{ i, j \}}$, denoted $X_i \indep X_j | \vect{X}_{[n] \setminus \{ i, j \}}$. In the particular case where $\vect{X}$ is a multivariate normal distribution, we refer to Gaussian graphical models \gcite{Uhler-2017}. Let then $\vect{X}$ be a Gaussian graphical model characterized by its covariance matrix $\matr{\Sigma} = ( \Sigma_{ij} )$, or, equivalently, its precision (or concentration) matrix $\matr{\Upsilon} = ( \Upsilon_{ij} ) = \matr{\Sigma} ^ { - 1 }$. Two other key quantities are the pairwise correlation matrix $\matr{\Omega} = ( \Omega_{ij} )$, defined as $\Omega_{ij} = \Sigma_{ij} / \sqrt{ \Sigma_{ii} \Sigma_{jj} }$ for $i \neq j$ and $\Omega_{ii} = 1$, as well as the partial correlation matrix $\matr{\Pi} = ( \Pi_{ij} )$, defined as $\Pi_{ij} = - \Upsilon_{ij} / \sqrt{ \Upsilon_{ii} \Upsilon_{jj} }$ for $i \neq j$ and $\Pi_{ii} = 1$. Then, for $i \neq j$, the relationship of conditional independence $X_i \indep X_j | \vect{X}_{[n] \setminus \{ i, j \}}$ is equivalent to $\Upsilon_{ij} = 0$ and $\Pi_{ij} = 0$ \gcite[Chap.~6]{Whittaker-1990}.
\par
Our interest in Gaussian graphical models originates from statistical mechanics, where the Ising model and its various extensions (Potts model, XY model, Heisenberg model, $n$-vector model, $\phi^4$ model) are used to investigate the behavior of variables related through various connection patterns. One extension of the Ising model to continuous real variables with noncompact support is the so-called Gaussian free field model \gcite[Chap.~8]{Friedli-2017}. In this case, each vertex $i \in \mathbb{Z}^d$ is associated with a real-valued variable $x_i$ and the corresponding Hamiltonian is of the form
$$\frac{\beta}{4d} \sum_{ i, j \in \mathbb{Z}^d : \| i - j \|_2 = 1 } ( x_i - x_j ) ^ 2 + \frac{ m ^ 2 } { 2 } \sum_{ i \in \mathbb{Z}^d } x_i ^ 2,$$
where $\beta \geq 0$ is the inverse temperature and $m \geq 0$ is the mass. In massive models ($m > 0$), pairwise correlation is known to decrease exponentially with distance \gcite[Prop.~8.30]{Friedli-2017}. 
\par
While this result is shown in the ``thermodynamic limit'', that is, for an infinite-dimensional variable (i.e., on $\mathbb{Z}^d$), we are here interested in the finite case. The reasons for this interest are twofold. First, a main way to approach statistical mechanics is through simulations, which only deal with finite case scenarios. It is therefore important to understand what the expected behavior of the system should be in such cases. Does pairwise correlation also decay exponentially? Also, we would like to gain a sense of how convergence from the finite to the infinite case occurs through some results regarding the speed of convergence.
\par
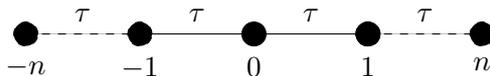
\begin{figure}[!htbp]
 \centering
 \begin{tikzpicture}[scale=1.5]
  \filldraw[dashed] (0,0) circle (3pt) node[below=5pt] {$-n$} -- node[above=1pt] {$\tau$} (1,0) circle (3pt) node[below=5pt] {$-1$};
  \filldraw (1,0) -- node[above=1pt] {$\tau$} (2,0) circle (3pt) node[below=5pt] {$0$} -- node[above=1pt] {$\tau$} (3,0) circle (3pt) node[below=5pt] {1};
  \filldraw[dashed] (3,0) circle (3pt) -- node[above=1pt] {$\tau$} (4,0) circle (3pt) node[below=5pt] {$n$};
 \end{tikzpicture}
 \caption{A conditional independence graph whose limit when $n \to \infty$ yields the one-dimensional Gaussian free field.} \label{fig:1}
\end{figure}

In the present study, we focus on the unidimensional case ($d = 1$) and consider the particular case of a (finite) Gaussian graphical model on $\vect{X}$ with an equicorrelational one-dimensional connection pattern between the $X_i$'s, as represented in Figure~\ref{fig:1}. Such a conditional independence graph entails that the Gaussian graphical model has a tridiagonal partial correlation matrix with an off-diagonal element $\tau$ that can be related to the parameters of the one-dimensional Gaussian free field by
\begin{equation} \label{eq:tau:gff}
 \tau = \frac{ \frac{ \beta } { 4 d } } { \frac{ 2 \beta } { 4 d } + \frac{ m ^ 2 } { 2 } }.
\end{equation}
We here restrict ourselves to the case $\tau > 0$ and only consider diagonally dominant matrices, leading to $0 \leq \tau < 1 / 2$ (which corresponds to the massive case, $m > 0$).

Under these assumptions, we show that $\Omega_{ij}^{(n)}$, the pairwise correlation between any two variables $X_i^{(n)}$ and $X_j^{(n)}$, decreases exponentially with the distance $| j - i |$ between variables, with a rate given by
\begin{equation} \label{eq:lambda}
 \lambda =  \arg \cosh \left( \frac{1} { 2 \tau } \right).
\end{equation}
More specifically, we show the following theorem.
\begin{theorem}
  Let $\vect{X}^{(n)}$ be a Gaussian graphical model with conditional independence graph given by Figure~\ref{fig:1}. Then the following results yield:
 \begin{itemize}
  \item $0 < \Omega_{ij}^{(n)} <  e ^ { - | j - i | \lambda }$ for all $n$;
  \item $\Omega_{ij}^{(n)} \to e ^ { - | j - i | \lambda }$ when $n \to \infty$;
  \item The absolute error $\Omega_{ij}^{(n)} - e ^ { - | j - i | \lambda }$ is $O \left[ e ^ { - 2 ( n + 1 ) \lambda } \right]$ when $n \to \infty$;
  \item The relative error $\Omega_{ij}^{(n)} e ^ { | j - i | \lambda } - 1$ is equal to
  $$- \left\{ \sinh [ 2 \max ( i, j ) \lambda ] -  \sinh [ 2 \min (i, j ) \lambda ] \right\} e ^ { - 2 ( n + 1 ) \lambda } + o \left[ e ^ { - 2 ( n + 1 ) \lambda } \right]$$
  when $n \to \infty$.
 \end{itemize}
\end{theorem}
Here, $O ( \cdot )$ and $o ( \cdot )$ are the usual big-O and little-o Bachmann--Landau notations, respectively, with
$$u_n = O ( v_n ) \qquad \Leftrightarrow \qquad \exists n_0, c \quad | u_n | < c | v_n | \ \forall n > n_0$$
and 
$$u_n = o( v_n ) \qquad \Leftrightarrow \qquad \frac{ u_n } { v_n } \stackrel{n \to \infty}{\to} 0.$$

\section{Proof of Theorem}

\begin{figure}[!htbp]
 \centering
 \begin{tikzpicture}[scale=1.5]
  \filldraw (0,0) circle (3pt) node[below=5pt] {1} -- node[above=1pt] {$\tau$} (1,0) circle (3pt) node[below=5pt] {2} -- node[above=1pt] {$\tau$} (2,0) circle (3pt) node[below=5pt] {3} -- node[above=1pt] {$\tau$} (3,0) circle (3pt) node[below=5pt] {4};
  \filldraw[dashed] (3,0) circle (3pt) -- node[above=1pt] {$\tau$} (4,0) circle (3pt) node[below=5pt] {$n$};
 \end{tikzpicture}
 \caption{Conditional independence graph of the Gaussian graphical model $\vect{Y}^{(n)}$.} \label{fig:2}
\end{figure}
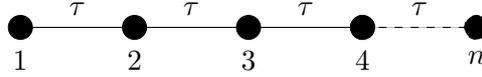

We start by expressing pairwise correlation in the case of the simpler model of an $n$-dimensional Gaussian graphical model $\vect{Y}^{(n)}$ with conditional independence graph given by Figure~\ref{fig:2}. We then relate the pairwise correlations for both models and derive the results for $\vect{X}^{(n)}$.

\subsection{Partial correlation matrix}

Assume that $\vect{Y}^{(n)}$ is a Gaussian graphical model with conditional independence graph given by Figure~\ref{fig:2}. The corresponding partial correlation matrix is then given by the following $n$-by-$n$ symmetric tridiagonal matrix
\begin{equation} \label{eq:corrpar}
 \Pi_{ij}^{(n)} = \left\{ \begin{array}{cl}
 1 & \mbox{if $i = j$} \\
 \tau & \mbox{if $| j - i | = 1$} \\
 0 & \mbox{otherwise.}
\end{array} \right.
\end{equation}

\subsection{From partial to pairwise correlation}

Letting $\matr{I}_n$ be the $n$-by-$n$ identity matrix and setting $\matr{\Upsilon}^{(n)} = 2 \matr{I}_n - \matr{\Pi}^{(n)}$, the (pairwise) correlation matrix $\matr{\Psi}^{(n)} = ( \Psi_{ij}^{(n)} )$ corresponding to the distribution can be obtained in two steps:
\begin{enumerate}
 \item Invert $\matr{\Upsilon}^{(n)}$ to obtain $\matr{\Sigma}^{(n)} = { \matr{\Upsilon}^{(n)} } ^ { - 1 }$;
 \item Decompose $\matr{\Sigma}^{(n)} = ( \Sigma_{ij}^{(n)} )$ using the correlation transform: $$\matr{\Sigma}^{(n)} = \matr{\Delta}^{(n)} \matr{\Psi}^{(n)} \matr{\Delta}^{(n)},$$
 where $\matr{\Delta} ^{(n)}= ( \Delta_{ij}^{(n)} )$ is a diagonal matrix with $\Delta_{ii}^{(n)} = \sqrt{\Sigma_{ii}^{(n)}}$.
\end{enumerate}

\subsection{Expression of $\Psi_{ij}^{(n)}$}

If $\matr{\Pi}^{(n)}$ has the form of Equation~(\ref{eq:corrpar}), then $\matr{\Upsilon}^{(n)}$ is also a tridiagonal matrix with off-diagonal element equal to $-\tau$. Defining $\lambda$ as in Equation~(\ref{eq:lambda}) and applying results from \gcitet{Hu_GY-1996}, we obtain that
$$\Sigma_{ij}^{(n)} = \frac{ 1 } { \tau } \, \frac{ \cosh [ ( n + 1 - | j - i | ) \lambda ] - \cosh [ ( n + 1 - i - j ) \lambda ] }{ 2 \sinh ( \lambda ) \, \sinh [ ( n + 1 ) \lambda ]}.$$
Using a basic identity of hyperbolic functions \gcite[ \S1.314]{Gradshteyn-2007}
$$\cosh ( x ) - \cosh ( y ) = 2 \sinh \left( \frac{x+y}{2} \right) \sinh \left( \frac{x-y}{2} \right),$$
we obtain
$$\Sigma_{ij}^{(n)} = \frac{ 1 } { \tau } \, \frac{ \sinh \left[ \frac{ 2 ( n + 1 ) - i - j - | j - i | } { 2 } \lambda \right] \, \sinh \left[ \frac{ i + j - | j - i | } { 2 } \lambda \right] }{ \sinh ( \lambda ) \, \sinh [ ( n + 1 ) \lambda ]}.$$
In particular, the diagonal elements read
$$\Sigma_{ii}^{(n)} = \frac{ 1 } { \tau } \, \frac{ \sinh \left[ ( n + 1 - i ) \lambda \right] \, \sinh \left( i \lambda \right) }{ \sinh ( \lambda ) \, \sinh [ ( n + 1 ) \lambda ]}.$$
This leads to the following expression for the correlation coefficient
$$\Psi_{ij}^{(n)} = \frac{ \sinh \left[ \frac{ 2 ( n + 1 ) - i - j - | j - i | } { 2 } \lambda \right] \, \sinh \left[ \frac{ i + j - | j - i | } { 2 } \lambda \right] }{ \sqrt{ \sinh \left[ ( n + 1 - i ) \lambda \right] \, \sinh \left( i \lambda \right) } \sqrt{ \sinh \left[ ( n + 1 - j ) \lambda \right] \, \sinh \left( j \lambda \right) } }.$$
In the following, we will restrict our attention to $i < j$ without loss of generality. For $j < i$, we can then use the symmetry identity $\Psi_{ij}^{(n)} = \Psi_{ji}^{(n)}$. So, if $j >i$, the previous result can be simplified to yield
\begin{eqnarray} \label{eq:Psi}
 \Psi_{ij}^{(n)} & = & \frac{ \sinh \left[ ( n + 1 - j ) \lambda \right] \, \sinh \left( i \lambda \right) }{ \sqrt{ \sinh \left[ ( n + 1 - i ) \lambda \right] \, \sinh \left( i \lambda \right) } \sqrt{ \sinh \left[ ( n + 1 - j ) \lambda \right] \, \sinh \left( j \lambda \right) } } \nonumber \\
 & = & \sqrt{ \frac{ \sinh \left[ ( n + 1 - j ) \lambda \right] \, \sinh \left( i \lambda \right) }{ \sinh \left[ ( n + 1 - i ) \lambda \right] \, \sinh \left( j \lambda \right) } } \nonumber \\
 & = & e ^ { - \lambda ( j - i ) } \sqrt{ \frac{ 1 -  e ^ { - 2 [ ( n + 1 - j ) ] \lambda } } { 1 -  e ^ { - 2 [ ( n + 1 - i ) ] \lambda } } \frac{ 1 -  e ^ { - 2 i \lambda } } { 1 -  e ^ { - 2 j \lambda } } }.
\end{eqnarray}

\subsection{Connection between $\vect{Y}^{(n)}$ and $\vect{X}^{(n)}$}

Gaussian free fields can be obtained as the limit when $n \to \infty$ of a $(2n+1)$-dimensional variables $\vect{X}^{(n)} = ( X_{-n}, \dots, X_{-1}, X_0, X_1, \dots, X_n )$ with a conditional independence graph given by Figure~\ref{fig:1}. Results regarding this model can be derived from the previous model and calculations by replacing $n$ with $2 n + 1$ and considering pairwise correlations of the form $\Omega_{ij}^{(n)} \equiv \Psi_{n+1+i,n+1+j}^{(2n+1)}$. In this perspective, Equation~(\ref{eq:Psi}) leads to, for $i < j$,
\begin{equation} \label{eq:Omega}
 \Omega_{ij}^{(n)} = e ^ { - ( j - i ) \lambda } \sqrt{ \frac{ 1 -  e ^ { - 2 [ ( n + 1 - j ) ] \lambda } } { 1 -  e ^ { - 2 [ ( n + 1 - i ) ] \lambda } } \frac{ 1 -  e ^ { - 2 ( n + 1 + i ) \lambda } } { 1 -  e ^ { - 2 ( n + 1 + j ) \lambda } } }.
\end{equation}

\subsection{Bounds}

From Equation~(\ref{eq:Omega}), it is straightforward to see that $\Omega_{ij}^{(n)}$ is always strictly positive. Also, since $u \mapsto 1 -  e ^ { - 2 ( n + 1 - u ) \lambda }$ is a strictly increasing function of $u$, and $u \mapsto 1 -  e ^ { - 2 ( n + 1 - u ) \lambda }$ a strictly decreasing function of $u$, we obtain for $i < j$
$$\sqrt{ \frac{ 1 -  e ^ { - 2 [ ( n + 1 - j ) ] \lambda } } { 1 -  e ^ { - 2 [ ( n + 1 - i ) ] \lambda } } } < 1 \qquad \mbox{and} \qquad \sqrt{ \frac{ 1 -  e ^ { - 2 ( n + 1 + i ) \lambda } } { 1 -  e ^ { - 2 ( n + 1 + j ) \lambda } } } < 1,$$
so that
$$0 < \Omega_{ij}^{(n)} < e ^ { -  ( j - i ) \lambda }$$
for all $n$.

\subsection{Asymptotics}

We can now provide the limit of $\Omega_{ij}^{(n)}$ when $n \to \infty$. Using the fact that $1 -  e ^ { - 2 [ ( n + 1 - u ) ] \lambda }$ tends to 1 when $n \to \infty$ for a given $u$, Equation~(\ref{eq:Omega}) leads to
\begin{equation} \label{eq:limite}
 \Omega_{ij}^{(n)} \stackrel{ n \to \infty } { \to } e ^ { - ( j - i ) \lambda }.
\end{equation}
Besides, using the following Taylor expansion for $u \to 0$,
\begin{equation} \label{eq:formasympt}
 ( 1 + u ) ^ k = 1 + k u + o ( u ),
\end{equation}
we can express $\Omega_{ij}^{(n)} / e ^ { - ( j - i ) \lambda }$ as
\begin{eqnarray*}
 \Omega_{ij}^{(n)} e ^ { ( j - i ) \lambda } & = & \left[ 1 -  e ^ { - 2 ( n + 1 - j ) \lambda } \right] ^ { \frac{1}{2} } \left[ 1 -  e ^ { - 2 ( n + 1 - i ) \lambda } \right] ^ { - \frac{1}{2} } \nonumber \\
 & & \qquad \times \left[ 1 -  e ^ { - 2 ( n + 1 + i ) \lambda } \right] ^ { \frac{1}{2} } \left[ 1 -  e ^ { - 2 ( n + 1 + j ) \lambda } \right] ^ { - \frac{1}{2} } \nonumber \\
 & = & \left\{ 1 - \frac{1}{2} e ^ { - 2 ( n + 1 - j ) \lambda } + o \left[ e ^ { - 2 ( n + 1 ) \lambda } \right] \right\} \nonumber \\
 & & \qquad \times \left\{ 1 + \frac{1}{2} e ^ { - 2 ( n + 1 - i ) \lambda } + o \left[ e ^ { - 2 ( n + 1 ) \lambda } \right] \right\} \nonumber \\
 & & \times \left\{ 1 - \frac{1}{2} e ^ { - 2 ( n + 1 + i ) \lambda } + o \left[ e ^ { - 2 ( n + 1 ) \lambda } \right] \right\} \nonumber \\
 & & \qquad \times \left\{ 1 + \frac{1}{2} e ^ { - 2 ( n + 1 + j ) \lambda } + o \left[ e ^ { - 2 ( n + 1 ) \lambda } \right] \right\} \nonumber \\
 & = & 1 - \frac{ 1 } { 2 } \left[ e ^ { 2 j \lambda } - e ^ { 2 i \lambda } +  e ^ { - 2 i \lambda } - e ^ { - 2 j \lambda } \right] e ^ { - 2 ( n + 1 ) \lambda } + o \left[ e ^ { - 2 ( n + 1 ) \lambda } \right] \\
 & = & 1 - \left[ \sinh ( 2 j \lambda ) -  \sinh ( 2 i \lambda ) \right] e ^ { - 2 ( n + 1 ) \lambda } + o \left[ e ^ { - 2 ( n + 1 ) \lambda } \right].
\end{eqnarray*}
We therefore have that 
$$\Omega_{ij}^{(n)}  e ^ { ( j - i ) \lambda } = 1 + O \left[ e ^ { - 2 ( n + 1 ) \lambda } \right],$$
so that
$$\Omega_{ij}^{(n)} - e ^ { - ( j - i ) \lambda } = e ^ { - ( j - i ) \lambda } \left[ \Omega_{ij}^{(n)} e ^ { ( j - i ) \lambda } - 1 \right] = O \left[ e ^ { - 2 ( n + 1 ) \lambda } \right].$$

\subsection{General results}

All results were proved for $i < j$. As mentioned earlier, the case $j < i$ can be solved by using the symmetry identity $\Omega_{ij}^{(n)} = \Omega_{ji}^{(n)}$. The most general results can therefore be expressed by replacing $i$ with $\min (i,j)$, $j$ with $\max(i,j)$, and $j-i$ with $|j-i|$, leading to
\begin{itemize}
 \item Bounds: $0 < \Omega_{ij}^{(n)} <  e ^ { - | j - i | \lambda }$ for all $n$;
 \item Limit: $\Omega_{ij}^{(n)} \to e ^ { - | j - i | \lambda }$ when $n \to \infty$;
 \item Asymptotic expansion: the absolute error $\Omega_{ij}^{(n)} - e ^ { - | j - i | \lambda }$ is $O \left[ e ^ { - 2 ( n + 1 ) \lambda } \right]$, and the relative error is given by
 \begin{eqnarray*}
  \Omega_{ij}^{(n)} e ^ { ( j - i ) \lambda } - 1 & = & - \left\{ \sinh [ 2 \max(i,j) \lambda ] -  \sinh [ 2 \min(i,j) \lambda ] \right\} e ^ { - 2 ( n + 1 ) \lambda } \\
  & & \qquad + o \left[ e ^ { - 2 ( n + 1 ) \lambda } \right].
 \end{eqnarray*}
\end{itemize}

\section{Discussion}

In the present manuscript, we considered a (finite-dimensional) Gaussian graphical model with the conditional independence graph depicted in Figure~\ref{fig:1}. We proved that the pairwise correlation decays exponentially at a rate given by $\lambda$ of Equation~(\ref{eq:lambda}). We also provided bounds for pairwise correlation as well as asymptotic expansions of the absolute and relative errors.

These results are in line with what is known about the one-dimensional Gaussian free field. Indeed, setting $\beta = 1$, pairwise correlation is known to be of the form $\exp ( - \xi_m | j - i | )$ with \gcite[Th.~8.33]{Friedli-2017}
$$\xi_m = \ln ( 1 + m^2 + \sqrt{ 2 m ^ 2 + m ^ 4 } ).$$
Using the relationship between $\tau$ and $( \beta, m )$ of Equation~(\ref{eq:tau:gff}) as well as the expression of $\arg \cosh$ in terms of logarithm \gcite[ \S1.622]{Gradshteyn-2007}, it can be shown that $\xi_m$ corresponds to our $\lambda$. 

Another quantity of interest is $\alpha = e^{-\lambda}$, which can be expressed using again the expression of $\arg \cosh$ in terms of logarithm \gcite[ \S1.622]{Gradshteyn-2007}, leading to
$$\frac{ 1 } { \alpha } = \frac{ 1 + \sqrt{ 1 - { 4 \tau ^ 2 } } } { 2 \tau },$$
or equivalently
\begin{equation} \label{eq:tridiag:alpha}
 \alpha = \frac{ 1 - \sqrt{ 1 - { 4 \tau ^ 2 } } } { 2 \tau }.
\end{equation}
From the definition, it is obvious that $\alpha \in [0,1)$, and that pairwise correlation decreases as $\alpha ^ { | j - i | }$. $\alpha$ appears naturally in the case where the Gaussian graphical model has a partial correlation matrix that is circulant instead of tridiagonal (see below).

One could wonder what the results are for the Gaussian graphical model $\vect{Y}^{(n)}$ with conditional independence graph of Figure~\ref{fig:2} that was used to derive our main results. The corresponding results are given in Appendix~\ref{an:1}. They are more complex due to the proximity of the boundary point 0 to $i$ and $j$.

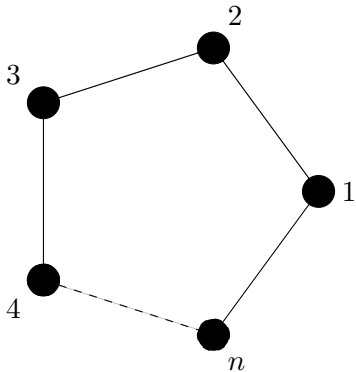
\begin{figure}[!htbp]
 \centering
 \begin{tikzpicture}[scale=2]
  \filldraw (1,0) circle (3pt) node[label={[label distance=1pt]0:1}] {} -- (360/5:1) circle (3pt) node[label={[label distance=1pt]360/5:2}] {} -- (360*2/5:1) circle (3pt) node[label={[label distance=1pt]360*2/5:3}] {} -- (360*3/5:1) circle (3pt) node[label={[label distance=1pt]360*3/5:4}] {};
  \filldraw[dashed] (360*3/5:1) circle (3pt) -- (360*4/5:1) circle (3pt) node[label={[label distance=1pt]360*4/5:$n$}] {};
  \draw (360*4/5:1) -- (1,0);
 \end{tikzpicture}
  \caption{Another instance of conditional independence graph, which corresponds to a symmetric circulant partial correlation matrix.} \label{fig:3}
\end{figure}

Another finite pattern of conditional independence that would lead to one-dimensional Gaussian free fields is the one given in Figure~\ref{fig:3}. In this case, the partial correlation matrix is symmetric circulant and it can be shown that the pairwise correlation still decays exponentially with the same rate $\lambda$ (see Appendix~\ref{an:2}). However, we were not able to provide bounds nor an asymptotic expansion in that particular case.

Our results show that pairwise correlation in (finite-dimensional) Gaussian graphical models behave in a manner very similar to one-dimensional (infinite-dimensional) Gaussian free fields, the difference between both cases decreasing exponentially with $n$. As a consequence, computer simulations can be trusted to provide precise approximations for the behavior of one-dimensional Gaussian free fields.

Beyond pairwise correlation, a measure that we think would be relevant to quantify the global level of dependence within the system is a multivariate generalization of mutual information known as total correlation \gcite{Watanabe-1960}, multivariate constraint \gcite{Garner-1962}, $\delta$ \gcite{Joe-1989}, or multiinformation \gcite{Studeny-1998}. In the case of multivariate normal distributions, this measure has a simple expression in terms of the covariance matrix. While we were not able to provide the closed form expression in the case of tridiagonal nor circulant partial correlation matrices, we believe that such expressions might be helpful to understand the global behavior of the system.

Now that we have solved the case $d = 1$, we would like to investigate more general cases with more complex connectivity patterns, still in the case of a finite $n$. Note that a major advantage of multivariate normal distributions is that their structures of conditional independence can be read off their precision matrices. For instance, moving from a one-dimensional to a two-dimensional model simply implies to change from a tridiagonal partial correlation matrix to a partial correlation matrix with more non-zero bands. More complex connectivity patterns with specific features (e.g., random or small world) simply translate into different patterns in the precision matrix which can then be investigated either analytically or through computer simulations. And, again, multiinformation has a simple expression that could provide interesting insight into the global behavior of the system.

\appendix

\section{Results for $\vect{Y}^{(n)}$} \label{an:1}

\subsection{Bounds}

We start from Equation~(4) of the manuscript. From this equation, it is obvious that $\Psi_{ij}^{(n)} > 0$. Since $u \mapsto 1 -  e ^ { - 2 u \lambda }$ is a strictly increasing function of $u$ and $u \mapsto 1 -  e ^ { - 2 [ ( n + 1 - u ) ] \lambda }$ a strictly decreasing function of $u$, the term in the square root is smaller than one for $i < j$ and
$$\Psi_{ij}^{(n)} <  e ^ { - ( j - i ) \lambda }.$$
We therefore still have that pairwise correlation decreases exponentially with distance.

\subsection{Limit}

We can now provide the limit of $\Psi_{ij}^{(n)}$ when $n \to \infty$. Still from Equation~(4) of the manuscript, we have
\begin{equation} \label{eq:tridiag:limite}
 \Psi_{ij}^{(n)} \stackrel{ n \to \infty } { \to } e ^ { - ( j - i ) \lambda } \, \sqrt{  \frac{ 1 -  e ^ { - 2 i \lambda } } { 1 -  e ^ { - 2 j \lambda } } } \equiv \Psi_{ij}^{(\infty)}.
\end{equation}
Note that, in this case, $\Psi_{ij}^{(\infty)}$ is a function of both $i$ and $j$ that cannot be expressed as a function of $j - i$ (the distance between $i$ and $j$) only. An upper bound for $\Psi_{ij}^{(\infty)}$ is given by
$$\Psi_{ij}^{(\infty)} < e ^ { - ( j - i ) \lambda }.$$
Also, still from Equation~(4) of the manuscript, we have
\begin{equation} \label{eq:tridiag:corr:2}
 \Psi_{ij}^{(n)} = \Psi_{ij}^{(\infty)} \sqrt{ \frac{ 1 -  e ^ { - 2 ( n + 1 - j ) \lambda } } { 1 -  e ^ { - 2 ( n + 1 - i ) \lambda }  } }.
\end{equation}
Since $u \mapsto 1 -  e ^ { - 2 ( n + 1 - u ) \lambda }$ is a strictly decreasing function of $u$, we obtain for $i < j$
$$\sqrt{ \frac{ 1 -  e ^ { - 2 ( n + 1 - j ) \lambda } } { 1 -  e ^ { - 2 ( n + 1 - i ) \lambda }  } } < 1,$$
so that
$$\Psi_{ij}^{(n)} < \Psi_{ij}^{(\infty)}$$
for all $n$.

\subsection{Asymptotic behavior}

Using the Taylor expansion of Equation~(7) of the manuscript, we can express $\Psi_{ij}^{(n)} / \Psi_{ij}^{(\infty)}$ as
\begin{eqnarray}
 \frac{ \Psi_{ij}^{(n)} } { \Psi_{ij}^{(\infty)} } & = & \left[ 1 -  e ^ { - 2 ( n + 1 - j ) \lambda } \right] ^ { \frac{1}{2} } \left[ 1 -  e ^ { - 2 ( n + 1 - i ) \lambda } \right] ^ { - \frac{1}{2} } \nonumber \\
 & = & \left\{ 1 - \frac{1}{2} e ^ { - 2 ( n + 1 - j ) \lambda } + o \left[ e ^ { - 2 ( n + 1 ) \lambda } \right] \right\} \nonumber \\
 & & \qquad \times \left\{ 1 + \frac{1}{2} e ^ { - 2 ( n + 1 - i ) \lambda } + o \left[ e ^ { - 2 ( n + 1 ) \lambda } \right] \right\} \nonumber \\
 & = & 1 - \frac{ e ^ { 2 j \lambda } - e ^ { 2 i \lambda } } { 2 } e ^ { - 2 ( n + 1 ) \lambda } + o \left[ e ^ { - 2 ( n + 1 ) \lambda } \right].
\end{eqnarray}
This result directly entails that
$$\Psi_{ij}^{(n)} - \Psi_{ij}^{(\infty)} = \Psi_{ij}^{(\infty)} \left[ \frac{ \Psi_{ij}^{(n)} } { \Psi_{ij}^{(\infty)} } - 1 \right] = O \left[ e ^ { - 2 ( n + 1 ) \lambda } \right].$$

\subsection{General results}

All previous results were proved for $i < j$. The case $j < i$ is solved by using the symmetry identity $\Psi_{ij}^{(n)} = \Psi_{ji}^{(n)}$, so that $i$, $j$ and $j-i$ are replaced with $\min(i,j)$, $\max(i,j)$ and $|j-i|$, respectively. In the end, setting
\begin{equation} \label{eq:tridiag:corr:lim}
 \Psi_{ij}^{(\infty)} = e ^ { - | j - i | \lambda } \, \sqrt{  \frac{ 1 -  e ^ { - 2 \min ( i, j ) \lambda } } { 1 -  e ^ { - 2 \max ( i, j ) \lambda } } },
\end{equation}
we obtain the following results:
\begin{itemize}
 \item Bounds: $0 < \Psi_{ij}^{(n)} < \Psi_{ij}^{(\infty)} <  e ^ { - | j - i | \lambda }$;
 \item Limit: $\Psi_{ij}^{(n)} \to \Psi_{ij}^{(\infty)}$ as $n \to \infty$;
 \item Asymptotic expansion:
 $$\frac{ \Psi_{ij}^{(n)} } { \Psi_{ij}^{(\infty)} } = 1 - \frac{1}{2} \left[ e ^ { 2 \max ( i, j ) \lambda } - e ^ { 2 \min ( i, j ) \lambda } \right] e ^ { - 2 ( n + 1 ) \lambda } + o \left[ e ^ { - 2 ( n + 1 ) \lambda } \right]$$
 and
 $$\Psi_{ij}^{(n)} - \Psi_{ij}^{(\infty)} = O \left[ e ^ { - 2 ( n + 1 ) \lambda } \right].$$
\end{itemize}

\section{Circulant partial correlation matrix} \label{an:2}

\subsection{Model}

Assume that $\vect{X}^{(n)}$ is a Gaussian graphical model with a conditional independence graph given by Figure~3 of the manuscript. The corresponding partial correlation matrix is then given by the following $n$-by-$n$ symmetric circulant matrix
\begin{equation} \label{eq:tri:corrpar}
 \Pi_{ij}^{(n)} = \left\{ \begin{array}{cl}
 1 & \mbox{if $i = j$} \\
 \tau & \mbox{if $| j - i | \in \{ 1, n - 1 \}$} \\
 0 & \mbox{otherwise.}
\end{array} \right.
\end{equation}
It can be expressed in the general form of circulant matrices as
$$\matr{\Pi}^{(n)} = \mathrm{circ} \left[ c_0^{(n)}, c_1^{(n)}, \cdots, c_{n-1}^{(n)} \right]$$
with $c_0^{(n)} = 1$, $c_1^{(n)} = c_{n-1}^{(n)} = \tau$ and 0 otherwise.

\subsection{Pairwise correlation}

In this case, $\matr{\Upsilon}^{(n)} = 2 \matr{I}_n - \matr{\Pi}^{(n)}$ is a symmetric circulant matrix as well with
$$\matr{\Upsilon}^{(n)} = \mathrm{circ} ( 1, -\tau, 0, \dots, 0, -\tau).$$
For $n \geq 2$, the $n$ eigenvalues of $\matr{\Upsilon}^{(n)}$ are given by \gcite{Chen_M-1987, Chao_CY-1988}
\begin{equation} \label{eq:circ:vp}
 \mu_k^{(n)} = 1 - 2 \tau \cos ( k \theta_n ), \qquad k = 0, \dots, n-1,
\end{equation}
where we set $\theta_n = 2 \pi /n$. Note that we have $\mu_0^{(n)} = 1 - 2 \tau$; for $k \geq 1$, $\mu_{n-k}^{(n)} = \mu_k^{(n)}$ ; for $n$ even, we also have $\mu_{ \frac{n}{2} }^{(n)} = 1 + 2 \tau$. Let $\matr{Q}^{(n)} = n ( \matr{\Upsilon}^{(n)} ) ^ { - 1 }$, so that $\matr{\Sigma}^{(n)} = ( \matr{\Upsilon}^{(n)} ) ^ { - 1 } = \frac{1}{n} \matr{Q}^{(n)}$. Then $\matr{Q}^{(n)}$ is also a symmetric circulant matrix,
$$\matr{Q}^{(n)} = \mathrm{circ} \left[ q_0^{(n)}, q_1^{(n)}, \dots, q_{ n - 1 }^{(n)} \right],$$
with \gcite{Chao_CY-1988}
\begin{equation} \label{eq:circ:qk}
 q_k^{(n)} = \sum_{ j = 0 } ^ { n - 1 } \frac{ e ^ { - i j k \theta_n } } { \mu_j^{(n)} } = \sum_{ j = 0 } ^ { n - 1 } \frac{ e ^ { - i j k \theta_n } } { 1 - 2 \tau \cos ( j \theta_n ) }.
\end{equation}
In particular, we have for the diagonal term ($k=0$)
\begin{equation} \label{eq:circ:q0}
 q_0^{(n)} = \sum_{ j = 0 } ^ { n - 1 } \frac{ 1 } { \mu_j^{(n)} } = \sum_{ j = 0 } ^ { n - 1 } \frac{ 1 } { 1 - 2 \tau \cos ( j \theta_n ) }.
\end{equation}
Since $\matr{Q}^{(n)}$ is a symmetric circulant matrix, so is $\matr{\Sigma}^{(n)}$,
$$\matr{\Sigma}^{(n)} = \mathrm{circ} \left[ \sigma_0^{(n)}, \sigma_1^{(n)}, \dots, \sigma_{ n - 1 }^{(n)} \right],$$
with 
\begin{equation} \label{eq:circ:cov}
  \sigma_k^{(n)} = \frac{ q_k^{(n)} } { n }.  
\end{equation}
Finally, the correlation matrix $\matr{\Omega}^{(n)}$ is also a symmetric circulant matrix,
$$\matr{\Omega} = ( \Omega_{ij}^{(n)} ) = \mathrm{circ} \left[ 1, \omega_1^{(n)}, \dots, \omega_{ n - 1 }^{(n)} \right]$$
with
\begin{equation} \label{eq:circ:corr}
 \omega_k^{(n)} = \frac{ \sigma_k^{(n)} } { \sqrt{ [ \sigma_0^{(n)} ] ^ 2 } } = \frac{ \sigma_k^{(n)} } { \sigma_0^{(n)} } = \frac{ q_k^{(n)} } { q_0^{(n)} },
\end{equation}
with $q_k^{(n)}$ and $q_0^{(n)}$ given by Equations (\ref{eq:circ:qk}) and (\ref{eq:circ:q0}), respectively.

\subsection{Riemannian sum}

Let $h_k$ be the function that maps any $x \in [ 0, 2 \pi ]$ to 
\begin{equation} \label{eq:circ:hk}
 h_k ( x ) = \frac{ e^{ - i k x } }{ 1 - 2 \tau \cos ( x ) }.
\end{equation}
Setting $x_j^{(n)} = j \theta_n$ for $j = 0, \dots, n$, we have
$$0 = x_0^{(n)} < x_1^{(n)} < \dots < x_n^{(n)} = 2 \pi.$$
We now define
\begin{equation} \label{eq:circ:Sk}
 S_k^{(n)} = \sum_{ j = 0 } ^ { n - 1 } h_k [ x_j^{(n)} ] [ x_{j+1}^{(n)} - x_j^{(n)} ] = \theta_n \sum_{ j = 0 } ^ { n - 1 } \frac{ e ^ { - i j k \theta_n } } { 1 - 2 \tau \cos ( j \theta_n ) } = \theta_n q_k^{(n)}.
\end{equation}
By construction, $S_k^{(n)}$ is a left Riemann sum that converges to
$$S_k^{(n)} \stackrel{n \to \infty}{\to} I_k = \int_0^{2\pi} h_k ( x ) \, \ud x.$$

\subsection{Computation of integral}

We therefore need to compute $I_k$. Using Euler's formula
$$e^{ix} = \cos ( x) + i \sin ( x ),$$
we obtain
$$I_k = \int_0^{2\pi} \frac{ e^{ - i k x } } { 1 - \tau \left( e^{ix} + e^{-ix} \right) } \, \ud x.$$
Performing the parameter change $z = e^{ix}$, we can now write this integral as a contour integral on the unit circle
\begin{eqnarray*}
 I_k & = & \oint_{ | z | = 1 } \frac{ \overline{z}^k } { 1 - \tau ( z + z ^ { - 1 } ) } \, \frac{ \ud z} { i z } \\
 & = & \frac{1}{i} \oint_{ | z | = 1 } \frac{ \overline{z}^k } { - \tau z ^ 2 + z - \tau } \, \ud z.
\end{eqnarray*}
The roots of $ - \tau z ^ 2 + z - \tau$ are given by
$$\alpha = \frac{ 1 - \sqrt{1 - 4 \tau ^ 2 } } { 2 \tau }$$
and $1/\alpha$. The integral therefore yields
$$I_k = - \frac{ 1 } { i \tau } \oint_{ | z | = 1 } \frac{ \overline{z}^k } { \left( z - \alpha \right) \left( z - \frac{1}{\alpha} \right) } \, \ud z.$$
Factoring the fraction yields
$$ \frac{ 1 } { \left( z - \alpha \right) \left( z - \frac{1}{\alpha} \right) } = \frac{ u } { z - \alpha} - \frac{ u } { z - \frac{1}{\alpha} }$$
with
$$u = \frac{ \alpha }{ \alpha^2 - 1 } = - \frac{ \tau } { \sqrt{ 1 - 4 \tau ^ 2 } }.$$
We therefore have for the integral
\begin{eqnarray*}
 I_k & = &  \frac{ 1 } { i \sqrt{ 1 - 4 \tau ^ 2 }  } \oint_{ | z | = 1 } \left( \frac{ \overline{z}^k } { z - \alpha } - \frac{ \overline{z} ^ k } { z - \frac{1}{\alpha} } \right) \, \ud z.
\end{eqnarray*}
$1/\alpha$ is outside the unit circle, so that
$$\oint_{ | z | = 1 } \frac{ \overline{z} ^ k } { z - \frac{1}{\alpha} } \, \ud z = 0.$$
For the other other integral, we need to compute the residual of $f ( z ) = \overline{z}^k / ( z - \alpha )$ at $z = \alpha$. Since it is a simple pole, we have 
$$\mathrm{Res}_{ z = \alpha } f ( z ) = \lim_{ z \to \alpha } ( z - \alpha ) f ( z ) = \lim_{ z \to \alpha } \overline{z} ^ k = \alpha ^ k,$$
since $\alpha \in \mathbb{R}$. We are then then led to
$$\oint_{ | z | = 1 } \frac{ \overline{z}^k } { z - \alpha } \, \ud z = 2 i \pi \alpha ^ k$$
and, finally,
\begin{equation} \label{eq:circ:Ik}
 I_k = \frac{ 2 \pi \alpha ^ k } { \sqrt{ 1 - 4 \tau ^ 2 } }.
\end{equation}
In particular, we have for $k = 0$
\begin{equation} \label{eq:circ:I0}
 I_0 = \frac{ 2 \pi } { \sqrt{ 1 - 4 \tau ^ 2 } }.
\end{equation}

\subsection{Asymptotic approximation}

Now that we computed $I_k$, we can go back to the pairwise correlation. Since the sum $S_k^{(n)}$ of Equation~(\ref{eq:circ:Sk}) converges toward the integral $I_k$, we have for $\sigma_k^{(n)}$, using Equations (\ref{eq:circ:cov}) and (\ref{eq:circ:Sk}),
\begin{equation}
 \sigma_k^{(n)} = \frac{ q_k^{(n)} } { n } = \frac{ S_k^{(n)} } { n \theta_n } = \frac{ S_k^{(n)} } { 2 \pi } \stackrel{n \to \infty}{\to}  \frac{ I_k } { 2 \pi } = \frac{ \alpha ^ k }{ \sqrt{ 1 - 4 \tau ^ 2 } }
\end{equation}
and for $\omega_k^{(n)}$, using Equations (\ref{eq:circ:corr}) and (\ref{eq:circ:Sk}),
\begin{equation}
 \omega_k^{(n)} = \frac{ q_k^{(n)} } { q_0^{(n)} } = \frac{ S_k^{(n)} } { S_0^{(n)} } \stackrel{n \to \infty}{\to} \frac{ I_k } { I_0 } = \alpha ^ k.
\end{equation}
Since $k = | j - i |$, we can conclude that
$$\Omega_{ij}^{(n)} \stackrel{n \to \infty}{\to} \alpha^{|j-i|}.$$


\end{document}